\documentclass{amsart}

\usepackage{mathtools,amsmath,amsthm,amsfonts,latexsym,amssymb,mathrsfs,setspace,verbatim,enumerate,cite,mathtools,tikz-cd,appendix,comment}
\usepackage{fullpage}

\mathtoolsset{showonlyrefs=true}

\newtheorem{thm}{Theorem}

\newtheorem{lem}[thm]{Lemma}

\numberwithin{thm}{section}
\numberwithin{equation}{section}

\newcommand{\rat}{\mathbb Q}
\newcommand{\real}{\mathbb R}

\newcommand{\alg}{\overline\rat}
\newcommand{\algt}{\alg^{\times}}

\newcommand{\intg}{\mathbb Z}
\newcommand{\nat}{\mathbb N}

\newcommand{\tors}{\mathrm{tors}}

\newcommand{\norm}{\mathrm{Norm}}

\newcommand{\gal}{\mathrm{Gal}}

\newcommand{\stab}{\mathrm{Stab}}

\newcommand{\comp}{\mathrm{com}}
\newcommand{\supp}{\mathrm{supp}}

\newcommand{\fin}{\mathrm{fin}}

\newcommand{\xx}{{\bf x}}

\begin{document}

\title[A classification of rational valued linear functionals]{A classification of $\mathbb Q$-valued linear functionals on $\overline{\mathbb Q}^\times$ modulo units}

\author[C.L. Samuels]{Charles L. Samuels}
\address{Christopher Newport University, Department of Mathematics, 1 Avenue of the Arts, Newport News, VA 23606, USA}
\email{charles.samuels@cnu.edu}
\subjclass[2020]{11R04, 46E10 (Primary), 11A25, 11G50, 11R32, 46E15 (Secondary)}
\keywords{Linear Functionals, Weil Height, Prime Omega Function, Galois Action}

\begin{abstract}
	Let $\overline{\mathbb Q}$ be an algebraic closure of $\rat$ and let $\overline \intg$ denote the ring of algebraic integers in $\overline{\mathbb Q}$.  If $\mathcal V = \overline{\mathbb Q}^\times/{\overline \intg}^\times$ then $\mathcal V$ 
	is a vector space over $\mathbb Q$.  We provide a complete classification all elements in the algebraic dual $\mathcal V^*$  of $\mathcal V$ in terms of another $\mathbb Q$-vector space called the space of consistent maps.
	With an appropriate norm on $\mathcal V$, we further classify the continuous elements of $\mathcal V^*$.  As applications of our results, we classify extensions of the prime Omega function 
	to $\mathcal V$ and discuss a natural action of the absolute Galois group $\gal(\overline{\mathbb Q}/\mathbb Q)$ on $\mathcal V$.  All results are closely connected to a 2009 article of Allcock and Vaaler on the
	vector space $\mathcal G = \overline{\mathbb Q}^\times/\overline{\mathbb Q}^\times_\tors$.
\end{abstract}

\maketitle

\section{Introduction}

\subsection{Background} Let $\alg$ be a fixed algebraic closure of $\rat$ and let $\algt_\tors$ denote the group of roots of unity in $\algt$.  Following the notation of \cite{AllcockVaaler}, we write $\mathcal G = \algt/\algt_\tors$
and note that $\mathcal G$ is a vector space over $\rat$ with addition given by
\begin{equation} \label{GAdd}
	\mathcal G \times \mathcal G \to \mathcal G,\quad (\alpha,\beta) \mapsto \alpha\beta
\end{equation}
and scalar multiplication given by
\begin{equation} \label{GMultiply}
	\rat\times \mathcal G \to \mathcal G,\quad (r,\alpha) \mapsto \alpha^r.
\end{equation}
For each number field $K$, we write $M_K$ to denote the set of all places of $K$.   If $v\in M_K$ we let $p_v$ be the unique place of $\rat$ such that $v$ divides $p_v$.
We write $K_v$ for the completion of $K$ with respect to $v$, and if $L$ is a finite extension of $K$ we have the well-known formula
\begin{equation} \label{DegreeSum}
	[L:K] = \sum_{w\mid v} [L_w:K_v],
\end{equation}
where $w$ runs over places of $L$ dividing $v$.  We let $\|\cdot \|_v$ be the unique extension to $K_v$ of the usual $p_v$-adic absolute value on $\rat_{p_v}$ and
define a second absolute value $|\cdot |_v$ on $K_v$ by
\begin{equation*}
	|x|_v = \|x\|_v^{[K_v:\rat_{p_v}]/[K:\rat]}.
\end{equation*}
Given a point $\gamma\in K$, we define
\begin{equation} \label{GNorm}
	\|\gamma\|_{\mathcal G} = \sum_{v\in M_K} \big| \log |\gamma|_v\big|.
\end{equation}
Because of our normalization of absolute values, the right hand side of \eqref{GNorm} does not depend on $K$, and moreover, its value is unchanged when $\gamma$ is multiplied by a root of unity.  
The basic properties of absolute values now imply that $\|\cdot\|_{\mathcal G}$ is a norm on $\mathcal G$ with respect to the usual absolute value on $\rat$.  It is worth noting that the product formula for $|\cdot |_v$ implies that
\begin{equation*}
	\frac{1}{2} \|\gamma\|_{\mathcal G} = \sum_{v\in M_K}\log^+|\gamma|_v.
\end{equation*}
The right hand side of this formula is commonly known as the {\it Weil height}.  

In view of our above observations, $\mathcal G$ has a completion with respect to $\|\cdot\|_{\mathcal G}$ that is a Banach space over $\real$.
A powerful result of Allcock and Vaaler \cite[Theorem 1]{AllcockVaaler} showed how to represent this completion as a certain subspace of an $L^1$ space.  Specifically, let $Y$ denote the set of all places of $\alg$, and for each number field
$K$ and each place $v$ of $K$ write $Y(K,v) = \left\{ y\in Y: y\mid v\right\}.$  For each $y\in Y(K,v)$ there exists an absolute value $\|\cdot\|_y$ on $\alg$ that extends the absolute value $\|\cdot\|_v$ on $K$.
The authors of \cite{AllcockVaaler} defined a locally compact Hausdorff topology on $Y$ having the set
$\{ Y(K,v): [K:\rat] < \infty,\ v\in M_K\}$ as a basis.  They further defined a regular measure $\lambda$ on the Borel sets $\mathcal B$ of $Y$ with the property that
\begin{equation*}
	\lambda(Y(K,v)) = \frac{[K_v:\rat_{p_v}]}{[K:\rat]}
\end{equation*}
for each place $v$ of a number field $K$.  For each $\gamma\in \mathcal G$, we let $f_\gamma:Y\to \real$ be given by $f_\gamma(y) = \log \|\gamma\|_y$.
Allcock and Vaaler proved \cite[Theorem 1]{AllcockVaaler} that $\gamma\mapsto f_\gamma$ is an isometric isomorphism of $\mathcal G$ onto a dense $\rat$-linear subspace $\mathcal F$ of
\begin{equation*}
	\mathcal X := \left\{ F\in L^1(Y,\mathcal B,\lambda): \int_Y F(y)d\lambda(y) = 0\right\}.
\end{equation*}Over the last several years, a variety of authors have applied techniques of functional analysis on $\mathcal X$ to establish new results about $\mathcal G$ 
(see \cite{AkhtariVaaler,FiliMinerNorms,FiliMinerDirichlet,GrizzardVaaler,Vaaler}, for example).



Let $Y^0$ be the set of non-Archimedean places of $\alg$ and define
\begin{equation*}
	\overline \intg = \left\{ \alpha\in \alg: \|\alpha\|_y \leq 1\mbox{ for all } y\in Y^0\right\}.
\end{equation*}
By the strong triangle inequality, $\overline \intg$ is a subring of $\alg$, and the group of units of $\overline \intg$ is equal to 
\begin{equation*}
	{\overline \intg}^\times =  \left\{ \alpha\in \algt: \|\alpha\|_y  = 1\mbox{ for all } y\in Y^0\right\}.
\end{equation*}
Clearly ${\overline \intg}^\times$ contains the torsion subgroup of $\algt$ and $\mathcal U := {\overline \intg}^\times/\algt_{\tors}$ is a subspace of $\mathcal G$.  
In the language of \cite{AllcockVaaler}, $\mathcal U$ is the set of points $\alpha\in \mathcal G$ such that $\supp(f_\alpha) \subseteq Y(\rat,\infty)$.
Now let $$\mathcal V = \algt/{\overline \intg}^\times,$$ so according to the third isomorphism theorem, there exists a group isomorphism from $\mathcal V$ to $\mathcal G/\mathcal U$.
Of course, $\mathcal V$ is also a vector space over $\rat$ under the operations analogous to \eqref{GAdd} and \eqref{GMultiply}.  We now easily verify that 
$$\mathcal V\simeq \mathcal G/\mathcal U$$ as vector spaces over $\rat$.  If $y$ is a non-Archimedean place of $\alg$,
then $\|\cdot\|_y$ is well-defined on $\mathcal V$.  For each $\alpha\in \mathcal V$, we may now define $f_\alpha:Y^0\to \real$ by $f_\alpha(y) = \log\|\alpha\|_y$.

As was the case with $\mathcal G$, there is a natural norm on $\mathcal V$.  For a number field $K$, we write $M^0_K$ to denote the set of non-Archimedean places of $K$.
If $\gamma\in K$ then we define
\begin{equation} \label{SNorm}
	\|\gamma\|_{\mathcal V} = \sum_{v\in M^0_K} \big| \log |\gamma|_v\big|.
\end{equation}
Again, because of our normalization of absolute values, the right hand side of \eqref{SNorm} does not depend on $K$.  Furthermore, its value remains the same when $\gamma$ is multiplied by an element
of ${\overline \intg}^\times$.  It now follows that $\|\cdot\|_{\mathcal V}$ is a norm on $\mathcal V$ with respect to the usual absolute value on $\rat$ and the completion of $\mathcal V$ with 
respect to $\|\cdot\|_{\mathcal V}$ is a Banach space over $\real$.  

By applying \cite[Theorem 1]{AllcockVaaler}, we can show that the map $\alpha\mapsto f_\alpha$ defines an isomorphism of $\mathcal V$ onto a dense $\rat$-linear subspace $\mathcal R$ of $L^1(Y^0, \mathcal B,\lambda)$.  
If $\Phi:L^1(Y^0, \mathcal B,\lambda)\to \real$ is a continuous linear map, then the Riesz Representation Theorem \cite[Theorem 6.16]{Rudin} yields an element $\phi\in L^\infty(Y^0,\mathcal B,\lambda)$ such that
\begin{equation} \label{ZRep}
	\Phi(F) = \int_{Y^0} F(y)\phi(y)d\lambda(y)
\end{equation}
for all $F\in L^1(Y^0, \mathcal B,\lambda)$.  Unfortunately, \eqref{ZRep} provides no property of $\phi$ that is equivalent to the assertion that $\Phi(\mathcal R)\subseteq \rat$.  
As such, it does not create an explicit representation theorem for the set $\mathcal R'$ of continuous $\rat$-linear functionals on $\mathcal R$.  Other versions of the Riesz Representation Theorem such 
as \cite[Theorem 2.14]{Rudin} face a similar problem and apply only to certain classes of linear functionals on $\mathcal R$.

The primary goal of this article is to provide that representation theorem for the entire algebraic dual $\mathcal V^*$ of $\mathcal V$.  To accomplish this task, it is necessary to create a new structure which 
we call the space of consistent maps and denote by $\mathcal J^*$.   Roughly speaking, $\mathcal J^*$ plays a similar role to that of $L^\infty (Y^0,\mathcal B,\lambda)$ in \eqref{ZRep}.   
Our main representation theorem, as well as a continuous version of it, are presented later in the introduction.  In Section \ref{Applications}, we discuss an application of our main theorems to the prime Omega function.
Moreover, we describe a Galois action on $\mathcal V^*$ and state three different properties of this action that are consequences of our main theorems.  The proofs of all results are provided in the subsequent sections.

\subsection{Main Results} For each number field $K$, we let $\mathcal O_K$ be the ring of integers in $K$ and let $\mathcal H_K$ denote the image of $K^\times$ in $\mathcal V$ under the canonical map.  
Clearly the kernel of this map equals $\mathcal O_K^\times$, and therefore, it defines an isomorphism $K^\times/\mathcal O_K^\times\to \mathcal H_K$. 
If $v$ is a non-Archimedean place of $K$, then the map $\alpha\mapsto \|\alpha\|_v$ is a well-defined group homomorphism from $\mathcal H_K$ to $\real^\times$.  
We let $\mathcal V_K$ be the subspace of $\mathcal V$ spanned by $\mathcal H_K$.  It is straightforward to verify that 
\begin{equation*}
	\mathcal V_K = \left\{ \alpha\in \mathcal V: \alpha^n\in \mathcal H_K\mbox{ for some } n\in \nat\right\},
\end{equation*} 
and if $L$ is a number field containing $K$ then $\mathcal V_K\subseteq \mathcal V_L$.  If $\alpha\in \mathcal V_K$ and $m,n\in\nat $ are such that $\alpha^m,\alpha^n\in \mathcal H_K$ then
$\|\alpha^m\|_v^n = \|\alpha^n\|_v^m.$  As a result, we may define an extension of $\|\cdot\|_v$ to $\mathcal V_K$ by
\begin{equation*}
	\|\alpha\|_v = \|\alpha^n\|_v^{1/n},
\end{equation*}
where $n$ is a positive integer such that $\alpha^n\in \mathcal H_K$. 
 
We now define the set
\begin{equation*}
	\mathcal J = \{(K,v): [K:\rat] < \infty,\ v\in M^0_K\}.
\end{equation*}
A map $c:\mathcal J\to \rat$ is called {\it consistent} if
\begin{equation} \label{Consistency}
	c(K,v) = \sum_{w\mid v} c(L,w)
\end{equation}
for all number fields $K$, all finite extensions $L/K$, and all non-Archimedean places $v$ of $K$.  If $c$ and $d$ are consistent maps and $r\in \rat$ then we define
\begin{equation*}
	(c + d)(K,v) = c(K,v) + d(K,v)\quad\mbox{and}\quad (rc)(K,v) = rc(K,v).
\end{equation*}
Under these operations, the set $\mathcal J^*$ of consistent maps is a vector space over $\rat$.  

The simplest example of a consistent map is the element $c\in \mathcal J^*$ given by
$c(K,v) = [K_v:\rat_{p_v}]/[K:\rat]$, in which case, definition \eqref{Consistency} reduces to the well-known formula \eqref{DegreeSum}.
In fact, $\mathcal J^*$ contains many additional familiar elements arising from \cite{AllcockVaaler}.  Indeed, let $\mathcal M^*$ be the $\rat$-vector space of signed Borel measures $\mu$ on $Y^0$ for which
$\mu(Y(K,v)) \in \rat$ for all $(K,v)\in \mathcal J$.
Then the elementary properties of measures ensure that $c(K,v) := \mu(Y(K,v))$ defines a unique consistent map.  Based on this observation, $\mathcal M^*$ may be interpreted a subspace of $\mathcal J^*$.

Using methods from the proof of Theorem \ref{ClosedSubgroup}, it is not difficult to show that $\mathcal M^*$ is a proper subset of $\mathcal J^*$.  Roughly speaking, a measure is required to satisfy {\bf countable} additivity for all
measurable subsets of $Y^0$ while consistent maps are only required to satisfy {\bf finite} additivity for the sets $Y(K,v)$.  We would be interested to see a classification of all consistent maps that do not arise from Borel measures.

While the relationship to measures is interesting, our primary motivation for definition \eqref{Consistency} is that consistency is precisely the condition needed to define a certain linear functional on $\mathcal V$.
For each element $c\in \mathcal J^*$, we define $\Phi_c:\mathcal V\to \rat$ by
\begin{equation} \label{PhiC}
	\Phi_c(\alpha) =  \sum_{v\in M^0_K} \frac{c(K,v)}{\log p_v}\log \|\alpha\|_v,
\end{equation}
where $K$ is any number field such that $\alpha\in \mathcal V_K$.
Because of the defining property of consistent maps \eqref{Consistency}, the right hand side of \eqref{PhiC} does not depend on $K$, and therefore, 
it defines a linear map from $\mathcal V$ to $\rat$.   As it turns out, all elements of $\mathcal V^*$ must have this form.

\begin{thm} \label{AllLinearClassification}
	The map $c\mapsto \Phi_c$ is a vector space isomorphism of $\mathcal J^*$ onto $\mathcal V^*$.
\end{thm}

As we noted in the introduction, we sought a representation theorem for all elements of the algebraic dual $\mathcal V^*$.  Theorem \ref{AllLinearClassification} provides precisely that result
in terms of consistent maps.  Fortunately, this theorem also provides a framework in which to classify elements of the continuous dual $\mathcal V'$.  To this end, we define
\begin{equation*}
	\mathcal J' = \left\{ c \in \mathcal J^*: \frac{[K:\rat]}{[K_v:\rat_{p_v}]}\cdot \frac{c(K,v)}{\log p_v}\mbox{ is bounded for } (K,v)\in \mathcal J\right\}.
\end{equation*}
It is straightforward to check that $\mathcal J'$ is a subspace of $\mathcal J^*$, and as it turns out, the elements of $\mathcal J'$ correspond precisely to the elements of $\mathcal V'$.

\begin{thm} \label{Continuous}
	The map $c\mapsto \Phi_c$ is a vector space isomorphism of $\mathcal J'$ onto $\mathcal V'$.
\end{thm}

If $\pi:\mathcal G\to \mathcal V$ denotes the canonical map, then $\pi$ is continuous with respect to the norms \eqref{GNorm} and \eqref{SNorm}.  Hence, the pullback $\pi^*$ of $\pi$
given by $\pi^*(\Phi) = \Phi\circ \pi$ is an injection of $\mathcal V^*$ into $\mathcal G^*$ and an injection of $\mathcal V'$ into $\mathcal G'$.  We would be very interested to see improved
versions of Theorem \ref{AllLinearClassification} and \ref{Continuous} that classify all elements of $\mathcal G^*$ and $\mathcal G'$, respectively.  As we shall find, our proofs cannot be easily
modified to reach this goal.

\section{Applications of our Representation Theorems} \label{Applications}

Let $\Omega:\rat^\times\to \intg$ be the unique group homomorphism such that $\Omega(p) = 1$ for every prime $p$.
This map, which we shall call the {\it classical prime Omega function}, has a very important role in many open problems in number theory.  
For example, Polya \cite{Polya} studied the function $L:\nat\to \intg$ given by
\begin{equation*} 
	L(x) = \sum_{n=1}^x (-1)^{\Omega(n)}.
\end{equation*}
Based on early numerical evidence, it seemed plausible that $L(x) \leq 0$ for all $x\geq 2$, and if this were correct, Polya showed that the Riemann 
Hypothesis would follow.  While this claim was later proven false \cite{Haselgrove,Lehman,Tanaka}, a strong connection between $L(x)$ and the Riemann 
Hypothesis remains intact.  Indeed, the Riemann Hypothesis is true if and only if $L(x) = O(x^{1/2 + \varepsilon})$ for every $\varepsilon > 0$ 
(see \cite{MossinghoffTrudgian} for the proof of an even stronger result).

In 1965, Chowla \cite{Chowla} studied the more general summation
\begin{equation*}
	L_f(x) = \sum_{n=1}^x (-1)^{\Omega(f(n))},
\end{equation*}
conjecturing that $L_f(x) = o(x)$ for all polynomials $f\in \intg[x]$ that do not have the form $c(g(x))^2$.  Chowla's conjecture follows from the Prime 
Number Theorem in the case that $f$ is linear, but all other cases remain open.  Cassaigne et al \cite{Cassaigne} posed the weaker conjecture that 
$L_f(x)$ changes sign infinitely many times, and subsequently, they established an assortment of special cases.   Borwein, Choi and Ganguli \cite{BCG} 
proved additional results involving certain quadratic polynomials.  Overall, progress toward Chowla's conjecture is extremely limited, and as such, 
it should be regarded as a very difficult problem even when $\deg f = 2$.

Clearly the classical prime Omega function is a well-defined map $\Omega: \mathcal H_\rat \to \intg$.   If $K$ is a number field and $\alpha\in K$, then we may define
\begin{equation} \label{CanonicalExtension}
	\Omega(\alpha) = \frac{\Omega(\norm_{K/\rat}(\alpha))}{[K:\rat]}.
\end{equation}
This definition is independent of $K$ and is unchanged when $\alpha$ is multiplied by an element of ${\overline \intg}^\times$.  Therefore, \eqref{CanonicalExtension} is a well-defined element of $\mathcal V^*$ 
that extends the classical prime Omega function.  Theorem \ref{AllLinearClassification} enables us to classify all elements of $\mathcal V^*$ that extend $\Omega$.

\begin{thm} \label{OmegaExtension}
	Suppose that $c:\mathcal J\to \rat$ is a consistent map.  Then $\Phi_c$ is an extension of $\Omega$ if and only if $c(\rat,p) = -1$ for all non-Archimedean places $p$ of $\rat$.
\end{thm}

We can now recognize the example \eqref{CanonicalExtension} in the context of Theorem \ref{OmegaExtension}.  We define the map $c:\mathcal J\to \rat$ by
\begin{equation} \label{BasicConsistent}
	c(K,v) = - \frac{[K_v:\rat_{p_v}]}{[K:\rat]}.
\end{equation}
By applying \eqref{DegreeSum}, we see clearly that $c$ is consistent, and moreover, we plainly have that $c(\rat,p) = -1$ for all non-Archimedean places $p$ of $\rat$.  The definition of $\Phi_c$ yields
\begin{equation*}
	\Phi_c(\alpha) = -\sum_{v\in M^0_K} \frac{\log |\alpha|_v}{\log p_v} = -\sum_{p\in M^0_\rat} \frac{1}{\log p} \sum_{v\mid p} \log |\alpha|_v
\end{equation*}
for any number field $K$ for which $\alpha\in \mathcal V_K$.  Because of the well-known identity
\begin{equation} \label{NormSwitch}
	\prod_{v\mid p} |\alpha|_v = \left| \norm_{K/\rat}(\alpha)\right|_p^{1/[K:\rat]}
\end{equation}
we conclude that
\begin{equation*}
	\Phi_c(\alpha) = -\sum_{p\in M^0_\rat} \frac{ \log| \norm_{K/\rat}(\alpha)|_p}{[K:\rat] \log p} =  \frac{\Omega(\norm_{K/\rat}(\alpha))}{[K:\rat]}.
\end{equation*}
In other words, selecting $c$ in accordance with \eqref{BasicConsistent} recreates the precise extension of $\Omega$ given by \eqref{CanonicalExtension}.  Interested readers will find an additional extension of $\Omega$ 
provided in Section \ref{GaloisApprox} immediately following the proof of Theorem \ref{WeakDense}.

Theorems \ref{AllLinearClassification} and \ref{Continuous} are useful tools to study an action of the absolute Galois group $G = \gal(\alg/\rat)$ on $\mathcal V^*$.
Before examining this action, we must briefly discuss an action of $G$ on $\mathcal V$.
If $\sigma\in G$ and $\alpha\in \mathcal V$ then $\sigma(\alpha)$ is a well-defined element of $\mathcal V$, and therefore, $(\sigma,\alpha) \mapsto \sigma(\alpha)$ is a left action of $G$ on $\mathcal V$.
We shall write
\begin{equation*}
	\stab_G(\alpha) = \left\{ \sigma\in G: \sigma(\alpha) = \alpha\right\}
\end{equation*}
for the stabilizer of $\alpha$ so that $\stab_G(\alpha)$ is certainly a subgroup of $G$.  When $G$ is equipped with the Krull topology \cite[Ch. IV, \S 1]{Neukirch} and $\mathcal V$ is given the discrete topology,
the action of $G$ on $\mathcal V$ is continuous, and therefore, $\stab_G(\alpha)$ is a closed subgroup of $G$.   For any closed subgroup $H$ of $G$, the definition of the Krull topology ensures that the following are equivalent:
\begin{enumerate}[(i)]
	\item\label{StabOpen} $H$ is an open subgroup of $G$.
	\item\label{StabFiniteIndex} $H$ has finite index in $G$.
	\item\label{StabNumberField} There exists a number field $K$ such that $\gal(\alg/K)\subseteq H$.
\end{enumerate}
It follows from \cite[Lemma 2.2]{FiliMinerNorms} or \cite[Lemma 3.2]{GrizzardVaaler} that these properties always hold when $H = \stab_G(\alpha)$.

The left action of $G$ on $\mathcal V$ induces a canonical right action of $G$ on $\mathcal V^*$ given by $(\Phi,\sigma)\mapsto \Phi\circ\sigma$. 
Now fix elements $\Phi\in \mathcal V^*$ and $\alpha\in \mathcal V$, and let $g: G\to \rat$ be given by $g(\sigma) = \Phi(\sigma(\alpha))$.  If $\mathcal V$ is given the discrete topology then $g$ is a composition of continuous maps
\begin{equation*}
	G\to G\times \mathcal V \to \mathcal V \to \rat.
\end{equation*}
Since we have
\begin{equation*}
	\stab_G(\Phi) = \bigcap_{\alpha\in \mathcal V}  \left\{ \sigma\in G: \Phi(\sigma(\alpha)) = \Phi(\alpha)\right\} = \bigcap_{\alpha\in \mathcal V} g^{-1}(\Phi(\alpha)),
\end{equation*}
we conclude that $\stab_G(\Phi)$ is a closed subgroup of $G$.  In contrast with the action of $G$ on $\mathcal V$, this subgroup need not have finite index in $G$. 

\begin{thm} \label{ClosedSubgroup}
	For each non-Archimedean place $p$ of $\rat$, let $x_p$ be a rational number.  Then there exists $c\in \mathcal J^*$ satisfying the following conditions:
	\begin{enumerate}[(i)]
		\item\label{RationalValues} $c(\rat,p) = x_p$ for all non-Archimedean places $p$ of $\rat$.
		\item\label{InfiniteIndex} $\stab_G(\Phi_c)$ has infinite index in $G$.
	\end{enumerate}
	Moreover, if $\{x_p/\log p: p\in M^0_\rat\}$ is bounded then there exists $c\in \mathcal J'$ satisfying \eqref{RationalValues} and \eqref{InfiniteIndex}.
\end{thm}

Of course, we may apply Theorem \ref{ClosedSubgroup} with $x_p = -1$ for all primes $p$.  Combined with Theorem \ref{OmegaExtension}, this special case yields an
extension of $\Omega$ to an element $\Phi\in \mathcal V'$ such that $\stab_G(\Phi)$ has infinite index in $G$.  

Although Theorem \ref{ClosedSubgroup} guarantees the existence of such elements, it is reasonable to ask whether they can be 
approximated by an element whose stabilizer has finite index.  To this end, we let
\begin{equation*}
	\mathcal V^*_{\fin} = \left\{ \Phi\in \mathcal V^*: [G:\stab_G(\Phi)] < \infty\right\}.
\end{equation*}
If $\Phi,\Psi\in \mathcal V^*$ and $r\in \rat$ then we certainly have that
\begin{equation*}
	\stab_G(\Phi) \cap \stab_G(\Psi) \subseteq \stab_G(\Phi + \Psi)\quad\mbox{and}\quad \stab_G(\Phi) \subseteq \stab_G(r\Phi),
\end{equation*}
and therefore, $\mathcal V^*_{\fin}$ is a subspace of $\mathcal V^*$. 

Now suppose $S$ is a set of non-Archimedean places of $\rat$.  We say that a consistent map $c:\mathcal J\to \rat$ is {\it supported on $S$} if $c(K,v) = 0$ for all number fields $K$ and all places
$v$ of $K$ not dividing a place in $S$.  Further, $c$ is called {\it compactly supported} if there exists a finite set $S$ such that $c$ is supported on $S$.  We say that $\Phi_c$ is {\it supported on $S$}
or {\it compactly supported} if $c$ satisfies the corresponding property.  We write $\mathcal V^*_{\comp}$ to denote the set of all compactly supported elements of $\mathcal V^*$ and note that $\mathcal V^*_{\comp}$ is a subspace
of $\mathcal V^*$.  We shall write $\mathcal C^* = \mathcal V^*_{\fin}\cap \mathcal V^*_{\comp}$.

\begin{thm} \label{WeakDense}
	$\mathcal C^*$ is a countable subspace of $\mathcal V'$.  Moreover, if $\Phi\in \mathcal V^*$ then there exists a sequence $\Phi_n\in \mathcal C^*$ which converges pointwise to 
	$\Phi$ with respect to the discrete topology on $\rat$.
\end{thm}

In addition to our efforts to approximate elements of $\mathcal V^*$ by elements of $\mathcal C^*$, Theorem \ref{ClosedSubgroup} raises another question: to determine whether every closed subgroup of $G$ is equal to 
$\stab_G(\Phi)$ for some $\Phi\in \mathcal V'$.  We have the following partial result. 

 \begin{thm} \label{OpenSubgroup}
 	Let $H$ be an open subgroup of $G$, and for each non-Archimedean place $p$ of $\rat$, let $x_p$ be a rational number.  Then there exists $c\in \mathcal J^*$ satisfying the following conditions:
	\begin{enumerate}[(i)]
		\item\label{RationalFixed} $c(\rat,p) = x_p$ for all non-Archimedean places $p$ of $\rat$
		\item\label{StabilizerFixed} $\stab_G(\Phi_c) = H$
	\end{enumerate}
	Moreover, if $\{x_p/\log p: p\in M^0_\rat\}$ is bounded then there exists $c\in \mathcal J'$ satisfying \eqref{RationalFixed} and \eqref{StabilizerFixed}
\end{thm}

\section{Proofs of our Representation Theorems} \label{RepTheoremProofs}

The proof of Theorem \ref{AllLinearClassification}, as well as several other results in this paper, rely on the existence a particular useful element in each number field $K$. 

\begin{lem} \label{SpecialPlace}
	Suppose that $K$ is a number field and $v$ is a non-Archimedean place of $K$.  Then there exists $\beta\in K$ such that $|\beta|_v < 1$ and $|\beta|_w = 1$ for all non-Archimedean places $w\ne v$.
	Moreover, if $v$ divides the rational prime $p\ne \infty$ then
	\begin{equation*}
		-e_v\log\|\beta\|_v = h\log p,
	\end{equation*}
	where $e_v$ is the ramification index of $v$ over $\rat$ and $h$ is the class number of $K$.
\end{lem}
\begin{proof}
 	Let $\mathcal O_K$ denote the ring of integers in $K$ and let 
	\begin{equation*}
		P = \{\alpha\in \mathcal O_K: |\alpha|_v < 1\}
	\end{equation*}
	so that $P$ is a prime ideal of $\mathcal O_K$.  We observe that $P^{h}$ is principal so we may let $\beta\in \mathcal O_K$ be a generator of $P^h$.  Clearly $\beta \in P$ so we immediately have that $|\beta|_v <1$.
	
	Now assume $w$ is a non-Archimedean place such that $|\beta|_w < 1$.  This means that $\beta$ belongs to the prime ideal
	\begin{equation*}
		Q := \{\alpha\in \mathcal O_K: |\alpha|_w < 1\},
	\end{equation*}
	and consequently,
	\begin{equation*}
		P^{h} = \beta\mathcal O_K \subseteq Q.
	\end{equation*}
	However, if $x\in P$ then $x^h\in P^h$, and hence, $x^h\in Q$.  But $Q$ is a prime ideal of $\mathcal O_K$, and therefore, $x\in Q$. We
	have shown that $P\subseteq Q$.  But every prime ideal in $\mathcal O_K$ is maximal, so this forces $Q = P$, and hence $w = v$.
	
	To prove the second statement, we apply the product formula to obtain
	\begin{equation*}
		|\beta|_v^{-1} = \prod_{w\mid\infty} |\beta|_w = \left|\norm_{K/\rat}(\beta)\right|_\infty^{1/[K:\rat]},
	\end{equation*}
	or equivalently,
	\begin{equation} \label{BetaToNorm}
		-[K_v:\rat_p]\log\|\beta\|_v = \log \left|\norm_{K/\rat}(\beta)\right|_\infty.
	\end{equation}
	Additionally, $\left|\norm_{K/\rat}(\beta)\right|_\infty$ is equal to the absolute norm of the ideal of $\mathcal O_K$ generated by $\beta$.  That is, if $N$ denotes the absolute norm on $K$, then 
	\begin{equation*}
		\left|\norm_{K/\rat}(\beta)\right|_\infty = N(P^h) = N(P)^h = p^{f_vh},
	\end{equation*}
	where $f_v = [\mathcal O_K/P: \mathbb F_p]$.  Combining this observation with \eqref{BetaToNorm} yields
	\begin{equation*}
		-[K_v:\rat_p]\log\|\beta\|_v = f_vh\log p
	\end{equation*}
	and the lemma follows from the fact that $[K_v:\rat_p] = e_vf_v$.
\end{proof}

Lemma \ref{SpecialPlace} will be used at various points throughout this article, but it plays a crucial role in our ensuing proof of Theorem \ref{AllLinearClassification}.

\begin{proof}[Proof of Theorem \ref{AllLinearClassification}]
	Suppose that $c,d\in \mathcal J^*$ and $r,s\in \rat$.  For each $\alpha\in \mathcal V$, we have that
	\begin{align*}
		\Phi_{rc + sd}(\alpha)  & = \sum_{v} \frac{rc(K,v) + sd(K,v)}{\log p_v} \log \|\alpha\|_v \\
			& = r \sum_{v} \frac{c(K,v)}{\log p_v} \log \|\alpha\|_v + s\sum_{v} \frac{d(K,v)}{\log p_v} \log \|\alpha\|_v \\
			& = (r\Phi_c + s\Phi_d)(\alpha),
	\end{align*}
	so it follows that $c\mapsto \Phi_c$ is linear.
	To see that this map is injective, we assume $c$ is a consistent map such that $\Phi_c\equiv 0$.  Given $(K,v)\in \mathcal J$, Lemma \ref{SpecialPlace} yields an element $\beta\in K$ such that
	$\|\beta\|_v < 1$ and $\|\beta\|_w = 1$ for all $w\in M^0_K \setminus\{v\}$.  Letting $\alpha$ be the image of $\beta$ in $\mathcal V$, we now obtain
	\begin{equation*}
		0 = \Phi_c(\alpha) = \frac{c(K,v)}{\log p_v} \log \|\beta\|_v
	\end{equation*}
	so that $c(K,v) = 0$.  As this property holds for all $(K,v)\in \mathcal J$, we deduce that $c\equiv 0$.
		
	In order to prove that our map is surjective, we define
	\begin{equation*}
		\mathcal T_K = \left\{ (a_v)_{v\in M^0_K}: a_v\in \rat,\ a_v=0\mbox{ for all but finitely many } v\in M^0_K\right\}
	\end{equation*}
	and note that $\mathcal T_K$ is a vector space over $\rat$.  For each $v\in M^0_K$, we define
	\begin{equation*}
		z_{v,w} = \begin{cases} 1 & \mbox{if } w = v \\ 0 &\mbox{if } w\ne v\end{cases}
	\end{equation*}
	and let
	\begin{equation*}
		{\bf z}_v = (z_{v,w})_{w\in M^0_K}.
	\end{equation*}
	We easily verify that $\{{\bf z}_v:v\in M^0_K\}$ is a basis for $\mathcal T_K$ over $\rat$, i.e., every element of $\mathcal T_K$ is expressed uniquely as a finite linear combination of the vectors
	${\bf z}_v$.
	
	Next, we let $\Delta_K:\mathcal V_K\to \mathcal T_K$ be given by
	\begin{equation*}
		\Delta_K(\alpha) = \left( \frac{\log \|\alpha\|_v}{\log p_v}\right)_{v\in M^0_K}.
	\end{equation*}
	By our definition of $\|\cdot\|_v$ on $\mathcal V_K$ given in the introduction, each component of $\Delta_K(\alpha)$ is indeed a rational number, and therefore, $\Delta_K$ is a well-defined map.
	We claim that $\Delta_K$ is an isomorphism.  Since each component of $\Delta_K$ defines a linear map, we obtain immediately that $\Delta_K$ is linear.  Now assume that $\alpha\in \mathcal V_K$ is such that 
	$\Delta_K(\alpha) = 0$.  There exists $n\in \nat$ such that $\alpha^n\in \mathcal H_K$, and by definition of $\Delta_K$, we obtain
	\begin{equation*}
		\|\alpha^n\|_v = 1
	\end{equation*}
	for all $v\in M^0_K$.  This means that $\alpha^n$ is the identity element in $\mathcal H_K$, and hence, $\alpha^n$ is the identity element in $\mathcal V_K$ as well.  Since $\mathcal V_K$ is a subspace of $\mathcal V$
	and $n\ne 0$, it follows that $\alpha$ is equal to the identity in $\mathcal V_K$.  We have now shown that $\Delta_K$ is injective.  To see that $\Delta_K$ is surjective, it is sufficient to prove that ${\bf z}_v$ is in the range of $\Delta_K$ 
	for all $v\in M^0_K$.  Let $\beta_v\in K$ be as in Lemma \ref{SpecialPlace} and let $\alpha_v$ be the image of $\beta_v$ in $\mathcal V_K$.
	We now have have
	\begin{equation*}
		\Delta_K(\alpha_v) = \frac{\log\|\beta_v\|_v}{\log p_v} {\bf z}_v,
	\end{equation*}
	and then applying the second statement of Lemma \ref{SpecialPlace}, we find that
	\begin{equation*}
		\Delta_K(\alpha_v) = -\frac{h}{e_v} {\bf z}_v,
	\end{equation*}
	where $h$ denotes the class number of $K$.  The surjectivity of $\Delta_K$ now follows from its linearity.
	
	Supposing that $(K,v)\in \mathcal J$, we let $L_K:\mathcal T_K\to \rat$ be given by $L_K = \Phi\circ \Delta_K^{-1}$ and define $c(K,v) = L_K({\bf z}_v)$. We claim that $\Phi_c = \Phi$.  
	Since $L_K$ is linear, we have
	\begin{equation} \label{LKBreakdown}
		L_K({\bf x}) = \sum_{v\in M^0_K} x_v L_K({\bf z}_v) = \sum_{v\in M^0_K} x_v c(K,v)
	\end{equation}
	for all ${\bf x} = (x_v)_{v\in M^0_K} \in \mathcal T_K$.  For an arbitrary point $\alpha\in \mathcal V$ we let $K$ be a number field such that $\alpha\in \mathcal V_K$.
	Now we apply \eqref{LKBreakdown} with ${\bf x} = \Delta_K(\alpha)$ to obtain
	\begin{equation} \label{FinalPhi}
		\Phi(\alpha) =  \sum_{v\in M^0_K} \frac{\log\|\alpha\|_v}{\log p_v}  c(K,v).
	\end{equation}
	As $\alpha$ was an arbitrary point in $\mathcal V$ and $K$ was an arbitrary number field with $\alpha\in \mathcal V_K$, this formula holds for all such $\alpha$ and $K$.
	
	In order to complete the proof of surjectivity, we need to show that $c$ is consistent.  To see this, let $L/K$ be an extension of number fields and let $v\in M^0_K$.  We once again apply Lemma \ref{SpecialPlace}
	to obtain an element $\beta\in K$ such that $\|\beta\|_v < 1$ and $\|\beta\|_w = 1$ for all $w\in M^0_K \setminus \{v\}$.  We may now use \eqref{FinalPhi} with $\alpha$ equal to the image of $\beta$ in $\mathcal V$
	and with both number fields $K$ and $L$.  We now find that
	\begin{equation*}
		 \frac{\log\|\beta\|_v}{\log p_v}  c(K,v) = \sum_{w\mid v} \frac{\log\|\beta\|_w}{\log p_w} c(L,w),
	\end{equation*}
	where the right hand sum runs over all places $w$ of $K$ dividing $v$.  Since $\|\beta\|_v = \|\beta\|_w$ and $p_v = p_w$, it immediately follows that
	\begin{equation*}
		c(K,v) = \sum_{w\mid v} c(L,w)
	\end{equation*}
	as is required to show that $c$ is consistent.
\end{proof}

We remarked following the statement of Theorem \ref{Continuous} that it would be difficult to adapt our methods to prove analogous representation theorems for $\mathcal G^*$ and $\mathcal G'$.
Now that we have seen the proof of Theorem \ref{AllLinearClassification}, it is clear why this is the case.  Indeed, twice during its proof, we required an element $\beta\in K$ such that $\|\beta\|_v \ne 1$
for precisely one place $v\in M^0_K$.  Certainly the product formula prohibits the existence of an element of $\beta \in K$ such that $\|\beta\|_v \ne 1$ for precisely one place $v\in M_K$.
Of course, it may still be possible that there exists a nontrivial point $\alpha\in \mathcal G$ such that $\supp(f_\alpha)\subseteq Y(K,v)$.  If such a point exists for every number field $K$ and every $v\in M_K$, 
then there is some hope of using our methods to classify all elements of $\mathcal G^*$.

The referee has pointed out an alternate proof of Theorem \ref{AllLinearClassification} which interprets $\mathcal V^*$ as inverse limit.  If $K\subseteq L$ are number fields then the natural inclusion
map $\mathcal V_K\hookrightarrow \mathcal V_L$ yields the direct limit
\begin{equation*}
	\mathcal V = \varinjlim_K \mathcal V_K.
\end{equation*}
As a result, the two dual spaces $\mathcal V^*$ and $\mathcal V'$ are inverse limits
\begin{equation*}
	\mathcal V^* = \varprojlim_K \mathcal V_K^*\quad\mbox{and}\quad \mathcal V' = \varprojlim_K \mathcal V_K'.
\end{equation*}
Using these observations along with \cite[Lemma 4]{AllcockVaaler}, the map $c$ may be built in a more straightforward way and its consistency follows from basic properties of inverse limits.
Although our earlier proof is somewhat more technical, it is more elementary.  Specifically, it does not rely on the deep results of \cite{AllcockVaaler} and does not require any knowledge of inverse limits.

\begin{proof}[Proof of Theorem \ref{Continuous}]
	It is sufficient to show that $\Phi_c$ is continuous if and only if there exists $C \geq 0$ such that
	\begin{equation} \label{StrongBound}
		\left| \frac{c(K,v)}{\log p_v}\right| \leq C\cdot \frac{[K_v:\rat_{p_v}]}{[K:\rat]}
	\end{equation}
	for all number fields $K$ and all non-Archimedean places $v$ of $K$.

	We first assume that \eqref{StrongBound} holds for all number fields $K$ and all non-Archimedean places $v$ of $K$.  For each $\alpha\in \alg$, we assume that $K$ is a number field containing $\alpha$
	so that
	\begin{equation*}
		|\Phi(\alpha)| \leq \sum_{v\in M^0_K} \left| \frac{c(K,v)}{\log p_v}\right|\cdot \big|\log \|\alpha\|_v \big| \leq C\sum_{v\in M^0_K} \frac{[K_v:\rat_p]}{[K:\rat]}\cdot \big|\log \|\alpha\|_v \big|
			= C\cdot \|\alpha\|_{\mathcal V}
	\end{equation*}
	and it follows that $\Phi$ is continuous with respect to $\|\cdot\|_{\mathcal V}$.
	
	Now suppose that $\Phi$ is continuous, and therefore, we may assume that $C\geq 0$ is such that $|\Phi(\alpha)| \leq C/2$ for all $\alpha$ belonging to the unit ball.
	We let $K$ be a number field and $v$ a place of $K$.  According to Lemma \ref{SpecialPlace}, there exists $\beta\in K$ such that $|\beta|_v <1$ and $|\beta|_w = 1$ for all places 
	$w\in M^0_K\setminus\{v\}$.  This means that
	\begin{equation*}
		\|\beta\|_{\mathcal V} =  \big| \log|\beta|_v\big| = \frac{[K_v:\rat_p]}{[K:\rat]} \cdot \big| \log\|\beta\|_v\big|.
	\end{equation*}
	Clearly $\|\beta\|_{\mathcal V}\ne 0$ so we may select a rational number $r$ such that
	\begin{equation} \label{rInequalities}
		\frac{1}{2} \leq r\cdot \|\beta\|_{\mathcal V} \leq 1.
	\end{equation}
	Since $\beta^r$ belongs to the unit ball in $\mathcal V$, we have that
	\begin{equation*}
		\frac{C}{2} \geq |\Phi(\beta^r)| = r\left| \frac{c(K,v)}{\log p_v}\right| \cdot \big| \log\|\beta\|_v\big|,
	\end{equation*}
	and then multiplying both sides by $[K_v:\rat_p]/[K:\rat]$ we obtain
	\begin{equation*}
		 \frac{[K_v:\rat_p]}{[K:\rat]} \cdot \frac{C}{2}  \geq r\cdot\left| \frac{c(K,v)}{\log p_v}\right| \big| \log|\beta|_v\big| = \left| \frac{c(K,v)}{\log p_v}\right|\cdot r\cdot \|\beta\|_{\mathcal V}.
	\end{equation*}
	The result now follows from \eqref{rInequalities}.
\end{proof}

The proof of Theorem \ref{OmegaExtension} is very brief, so we provide here before proceeding with the more technical proofs of Theorems \ref{ClosedSubgroup}, \ref{WeakDense} and \ref{OpenSubgroup}.

\begin{proof}[Proof of Theorem \ref{OmegaExtension}]
	We first assume that $\Phi(r) = \Omega(r)$ for all $r\in \mathcal H_\rat$.  For each rational prime $p\ne \infty$ we have that
	\begin{equation*}
		1 = \Phi(p) =  \frac{c(\rat,p)}{\log p}\log \|p\|_p = -c(\rat,p)
	\end{equation*}
	and we obtain the required property.  Now assuming that $c(\rat,p) = -1$ we let $r\in \rat^\times$ and write
	\begin{equation*}
		r = \prod_{i=1}^n p_i^{a_i}
	\end{equation*}
	for non-zero integers $a_i$.  Then we have
	\begin{equation*}
		\Phi(r) = -\sum_{p\in M^0_\rat} \frac{\log \|r\|_p}{\log p} = \sum_{i=1}^n a_i = \Omega(r)
	\end{equation*}
	as required.
\end{proof}

\section{Proof of Theorem \ref{ClosedSubgroup}} \label{ClosedSubgroupProof}

The proof of Theorem \ref{ClosedSubgroup} requires us to detect whether $\gal(\alg/K)\subseteq \stab_G(\Phi)$ for a given number field $K$.  Our first lemma
provides the necessary tools to do this when $\Phi$ is given in the form of Theorem \ref{AllLinearClassification}.
If $K$ is any field with $\rat\subseteq K\subseteq \alg$ and $\Phi\in \mathcal V^*$, we say that $\Phi$ is {\it $K$-Galois invariant} if $\gal(\alg/K)\subseteq \stab_G(\Phi)$.
 A consistent map $c:\mathcal J\to \rat$ is called {\it $K$-Galois invariant} if $\Phi_c$ is $K$-Galois invariant.  
If $K$ and $L$ are number fields with $K\subseteq L$, we follow the notation of \cite{AllcockVaaler} and write $W_v(L/K)$ to denote the set of places of $L$ that divide $v$.  If $L/K$ is a Galois extension and $w\in W_v(L/K)$, 
then we define the element $\sigma(w)\in W_v(L/K)$ by
\begin{equation*}
	\|\alpha\|_{\sigma(w)} = \|\sigma^{-1}(\alpha)\|_w.
\end{equation*}
The map $(\sigma,w)\mapsto \sigma(w)$ is transitive left action of $\gal(L/K)$ on $W_v(L/K)$.  

Before we provide our first lemma, we note several notational conventions that we will adopt throughout the rest of this article.  We will often use the symbols $K$ and $L$ to denote number fields.  Unless otherwise stated shall always assume, 
$p$ runs over places of $\rat$, $v$ runs over places of $K$, and $w$ runs over places of $L$.  In all cases, these symbols shall denote non-Archimedean places only.
If $\Phi\in \mathcal V^*$ and $\alpha\in \algt$, we shall simply write $\Phi(\alpha) = \Phi(\pi(\alpha))$, where $\pi:\algt\to \mathcal V$ is the canonical homomorphism.  If $c$ is a consistent map and $\alpha$ belongs to a number field $K$, then we certainly have that
\begin{equation*}
	\Phi_c(\alpha) = \sum_v \frac{c(K,v)}{\log p_v} \log \|\alpha\|_v.
\end{equation*}
In other words, the definition \eqref{PhiC} is unchanged when the element $\alpha\in\mathcal V$ is replaced with an element of $\algt$.

\begin{lem}\label{KGaloisInvariant}
	Suppose that $K$ is a number field and $c:\mathcal J\to \rat$ is a consistent map.  Then $c$ is $K$-Galois invariant if and only if
	\begin{equation} \label{cConstant}
		c(L,w) = \frac{[L_w:K_v]}{[L:K]} c(K,v)
	\end{equation}
	for all non-Archimedean places $v$ of $K$, all finite extensions $L/K$, and all places $w\in W_v(L/K)$.
\end{lem}
\begin{proof}
	We begin by assuming that \eqref{cConstant} holds for all non-Archimedean places $v$ of $K$, all finite extensions $L/K$, and all places $w\in W_v(L/K)$.  To prove that $\Phi$ is
	$K$-Galois invariant, let $\alpha\in \algt$ and $\sigma\in \gal(\alg/K)$.  We let $L$ be a finite Galois extension of $K$ containing $\alpha$ and note that the restriction of $\sigma$ to $L$ belongs
	to $\gal(L/K)$.  In our ensuing computations, we shall simply write $\sigma$ for this restriction.
	
	By our assumption that $L/K$ is Galois, the right hand side of \eqref{cConstant} depends only on $K, L$ and $v$.  We shall write $d(K,L,v)$ for this value so that
	\begin{equation*}
		d(K,L,v) = c(L,w)
	\end{equation*}
	for all places $w$ of $L$ dividing $v$, and therefore
	\begin{equation*}
		\Phi(\sigma(\alpha)) = \sum_{v}\sum_{w\mid v} \frac{c(L,w)}{\log p_v} \log\|\sigma(\alpha)\|_w = \sum_{v} \frac{d(K,L,v)}{\log p_v} \sum_{w\mid v} \log\|\alpha\|_{\sigma^{-1}(w)}.
	\end{equation*}
	The map $w\mapsto\sigma^{-1}(w)$ is a permutation of $W_v(L/K)$ so we obtain that
	\begin{equation*}
		\Phi(\sigma(\alpha)) =  \sum_{v} \frac{d(K,L,v)}{\log p_v} \sum_{w\mid v} \log\|\alpha\|_w = \sum_{v} \sum_{w\mid v} \frac{c(L,w)}{\log p_v} \log\|\alpha\|_w
			= \Phi(\alpha).
	\end{equation*}
	
	Now we assume that $\Phi$ is $K$-Galois invariant and that $L/K$ is a finite extension.  Assume that $E$ is a finite Galois extension of $K$ containing $L$, let $G = \gal(E/K)$, and let $v$ be a non-Archimedean
	place of $K$.  We shall adopt the convention that $t$ runs over non-Archimedean places of $E$.  We claim that
	\begin{equation} \label{ConstantOnV}
		c(E,t_1) = c(E,t_2) \quad\mbox{for all } t_1,t_2\in W_v(E/K).
	\end{equation}
	To see this, we apply Lemma \ref{SpecialPlace} to obtain an element $\beta\in E$ such that $\|\beta\|_{t_1} < 1$ and $\|\beta\|_t = 1$ for all non-Archimedean places $t\ne t_1$.
	Since $G$ acts transitively on $W_v(E/K)$, we may assume that $\sigma\in G$ is such that $\sigma(t_1) = t_2$.  These observations yield $\|\sigma(\beta)\|_{t_2} <1$ and
	$\|\sigma(\beta)\|_t = 1$ for all non-Archimedean places $t\ne t_2$.  Therefore, we have that
	\begin{equation*}
		\Phi(\beta) = \frac{c(E,t_1)}{\log p_v}\log\|\beta\|_{t_1}
	\end{equation*}
	and
	\begin{equation*}
		\Phi(\sigma(\beta)) = \frac{c(E,t_2)}{\log p_v} \log \|\sigma(\beta)\|_{t_2} = \frac{c(E,t_2)}{\log p_v} \log\|\beta\|_{t_1}
	\end{equation*}
	and \eqref{ConstantOnV} follows from $K$-Galois-invariance.
	
	We now know that $c(E,t)$ depends only on the place $v$ for which $t\mid v$.  Moreover, since $E/K$ is Galois, $[E_t:K_v]$ depends only on $v$ and $\#W_v(E/K) = [E:K]/[E_t:K_v]$.  Now applying the fact that
	$c$ is consistent, we obtain
	\begin{equation*}
		c(K,v) = \sum_{t\mid v} c(E,t) = c(E,t)\cdot \#W_v(E/K) = c(E,t)\cdot \frac{[E:K]}{[E_t:K_v]}
	\end{equation*}
	for all places $t\mid v$.  If $w$ is a place of $L$ dividing $v$, then we may assume that $t$ divides $w$.  By applying a similar argument we obtain
	\begin{equation*}
		c(L,w) = c(E,t) \frac{[E:L]}{[E_t:L_w]}
	\end{equation*}
	and the lemma follows.
\end{proof}

If $c$ is known to be $K$-Galois invariant and $L/K$ is a finite extension, then the values $c(L,w)$ are completely 
determined by \eqref{cConstant}.  Moreover, every other number field $F$ must be contained in a finite extension of $K$, and therefore, the values $c(F,u)$ are determined by the consistency of $c$.
In other words, a $K$-Galois invariant linear functional on $\mathcal V$ is completely determined by the values of $c(K,v)$. 

Among other things, the proof of Theorem \ref{ClosedSubgroup} requires us to find an element of $\mathcal V'$ which is not $K$-Galois invariant for any number field $K$.
This process begins with the following lemma concerning towers of number fields.  

\begin{lem} \label{GeneralExampleConstructor}
	Suppose that $\{K_i:i\in \nat\}$ is a collection of number fields such that
	\begin{equation*}
		\alg = \bigcup_{i=1}^\infty K_i\quad\mbox{and}\quad K_i\subseteq K_{i+1}
	\end{equation*}
	for all $i\in \nat$.  For each $i\in \nat$ and each non-Archimedean place $v$ of $K_i$ let $d_i(v)\in \rat$.  Finally, assume that
	\begin{equation} \label{dAscend}
		d_i(v) = \sum_{w\in W_v(K_{i+1}/K_i)} d_{i+1}(w)
	\end{equation}
	for all $i\in \nat$ and all non-Archimedean places $v$ of $K_i$.  Then there exists a unique consistent map $c:\mathcal J\to \rat$ such that $c(K_i,v) = d_i(v)$ for all $i\in \nat$ and all
	non-Archimedean places $v$ of $K_i$.
\end{lem}
\begin{proof}
	Let $F$ be a number field and let $u$ be a non-Archimedean place of $F$.  Since $F$ necessarily has the form $F = \rat(\gamma)$ for some $\gamma\in \alg$, there exists $i\in \nat$ such 
	that $F\subseteq K_i$.  Now let
	\begin{equation} \label{ciDef}
		c_i(F,u) = \sum_{v\in W_u(K_i/F)} d_i(v).
	\end{equation}
	It is straight forward to prove by induction using \eqref{dAscend} that
	\begin{equation*} \label{dJump}
		d_i(v) = \sum_{w\in W_v(K_j/K_i)} d_{j}(w)
	\end{equation*}
	for all $j\geq i$ and all non-Archimedean places $v$ of $K_i$.  Consequently, for each $j\geq i$ we have
	\begin{equation*}
		c_j(F,u) = \sum_{w\in W_u(K_j/F)} d_j(w) = \sum_{v\in W_u(K_i/F)} \sum_{w\in W_v(K_j/K_i)} d_j(w) = \sum_{v\in W_u(K_i/F)} d_i(v) = c_i(F,u).
	\end{equation*}
	In other words, the definition \eqref{ciDef} does not depend on $i$ and we shall simply write $c(F,u) = c_i(F,u)$.  To complete the proof, we must establish the following three claims:
	\begin{enumerate}[(A)]
		\item\label{Agreement2} $c(K_i,v) = d_i(v)$ for all $i\in \nat$ and all non-Archimedean places $v$ of $K_i$.
		\item\label{Consistent2} $c$ is consistent.
		\item\label{Unique2} $c:\mathcal J\to \rat$ is the unique map satisfying \eqref{Agreement2} and \eqref{Consistent2}.
	\end{enumerate}
	We obtain \eqref{Agreement2} immediately by applying \eqref{ciDef} with $F = K_i$ and $u = v$.  To establish \eqref{Consistent2}, we assume that $E$ is a finite extension of the number
	field $F$ and assume that $i\in \nat$ is such that $E\subseteq K_i$.  Then definition \eqref{ciDef} yields
	\begin{equation*}
		c(F,u) =  \sum_{v\in W_u(K_i/F)} d_i(v) =  \sum_{t\in W_u(E/F)} \sum_{v\in W_t(K_i/E)} d_i(v) = \sum_{t\in W_u(E/F)}  c(E,t).
	\end{equation*}
	Finally, assume that $c':\mathcal J\to \rat$ such that $c'(K_i,v) = d_i(v)$ for all $i\in \nat$ and all non-Archimedean places $v$ of $K_i$.  Let $F$ be a number field, let $u$ be a 
	non-Archimedean place of $F$, and let $i\in \nat$ be such that $F\subseteq K_i$.  Then we have
	\begin{equation*}
		c(F,u) =  \sum_{v\in W_u(K_i/F)} d_i(v) =  \sum_{v\in W_u(K_i/F)} c'(K_i,v) = c'(F,u),
	\end{equation*}
	where the last equality follows from our assumption that $c'$ is consistent.
\end{proof}

In our proof of Theorem \ref{ClosedSubgroup}, we will select values of $d_i(v)$ for use in Lemma \ref{GeneralExampleConstructor}.  As we shall find, it is necessary to avoid a situation in which 
$d_{i+1}(w) = d_i(v)$ for all $w$ dividing $v$.  Our next lemma shows that the fields $K_i$ may be selected in such a way.

\begin{lem} \label{NumberFieldChain}
	Suppose that $\{K_i:i\in \nat\}$ is a collection of number fields that are Galois over $\rat$ such that
	\begin{equation*}
		\alg = \bigcup_{i=1}^\infty K_i\quad\mbox{and}\quad K_i\subseteq K_{i+1}
	\end{equation*}
	for all $i\in \nat$.  If $v$ and $w$ are places of $K_i$ dividing the place $p$ of $\rat$ then $\#W_v(K_{i+1}/K_i) = \#W_w(K_{i+1}/K_i)$.  Moreover, there exist infinitely many
	values of $i$ such that $\#W_v(K_{i+1}/K_i) > 1$.
\end{lem}
\begin{proof}
	Suppose that $s$ and $t$ are places of $K_{i+1}$ dividing $v$ and $w$, respectively.  Since $K_{i+1}/K_i$ is Galois, we have that
	\begin{equation} \label{LocalSwitch}
		\#W_v(K_{i+1}/K_i) = \frac{[K_{i+1}:K_i]}{[(K_{i+1})_s:(K_i)_v]}\quad\mbox{and}\quad \#W_w(K_{i+1}/K_i) = \frac{[K_{i+1}:K_i]}{[(K_{i+1})_t:(K_i)_w]}.
	\end{equation}
	Similarly, since $K_{i+1}/\rat$ is Galois we obtain that $[(K_{i+1})_s:\rat_p] = [(K_{i+1})_t:\rat_p]$, and thus, we find that
	\begin{equation*}
		[(K_{i+1})_s:(K_i)_v] \cdot [(K_i)_v:\rat_p] = [(K_{i+1})_t:(K_i)_w] \cdot [(K_i)_w:\rat_p].
	\end{equation*}
	Using the fact that $K_i/\rat$ is Galois, we see that $[(K_i)_v:\rat_p] = [(K_i)_w:\rat_p]$ which means that $[(K_{i+1})_s:(K_i)_v] = [(K_{i+1})_t:(K_i)_w]$ so the first statement of the lemma
	follows from \eqref{LocalSwitch}.

	To prove the second statement, let $v_i$ be a place of $K_i$ such that $v_{i+1}\mid v_i$ for all $i \in \nat$.  In view of the first statement, it is sufficient to prove that there exist infinitely many
	values of $i$ such that $\#W_{v_i}(K_{i+1}/K_i) > 1$.  Supposing otherwise, there exists $n\in \nat$ such that $\#W_{v_i}(K_{i+1}/K_i) = 1$ for all $i\geq n$.
	Since $K_{i+1}/K_i$ is Galois, we have 
	\begin{equation*}
		\#W_{v_i}(K_{i+1}/K_i) = \frac{[K_{i+1}:K_i]}{[(K_{i+1})_{v_{i+1}}: (K_i)_{v_i}]}
	\end{equation*}
	so that $[(K_{i+1})_{v_{i+1}}: (K_i)_{v_i}] = [K_{i+1}:K_i]$ for all $i\geq n$.  Letting $\delta_i = [(K_i)_{v_i}:(K_1)_{v_1}]$ and $\Delta_i = [K_{i+1}:K_1]$ we obtain immediately that
	\begin{equation} \label{StableDelta}
		\frac{\delta_i}{\Delta_i} = \frac{\delta_n}{\Delta_n}
	\end{equation}
	 for all $i\geq n$.  Since $K_i/K_1$ is Galois, we know that $\delta_i\mid \Delta_i$ and $\Delta_i/\delta_i = \#W_{v_1}(K_i/K_1)$.  We know that $\#W_{v_1}(K_i/K_1)$ tends to $\infty$ as $i\to\infty$
	 (see \cite[Lemma 7.1]{KellySamuels} for a stronger result), and therefore
	 \begin{equation*}
	 	\lim_{i\to\infty} \frac{\delta_i}{\Delta_i}  = 0 
	\end{equation*}
	contradicting \eqref{StableDelta}.
\end{proof}

Equipped with our three crucial lemmas, we may proceed with our proof of Theorem \ref{ClosedSubgroup}.

\begin{proof}[Proof of Theorem \ref{ClosedSubgroup}]
	Since the set of all Galois extensions of $\rat$ is a countable set, we mat let $\{L_i:i\in \nat\}$ be the set of all such fields and assume that $L_1 = \rat$.  Now inductively define number fields $K_i$ by
	\begin{equation*}
		K_1 = L_1\quad\mbox{and}\quad K_{i+1} = K_{i}L_{i+1}\mbox{ for all } i\in \nat.
	\end{equation*}
	Certainly we have $K_1=\rat$, and additionally, 
	\begin{equation} \label{FullContain}
		\alg = \bigcup_{i=1}^\infty K_i\quad\mbox{and}\quad K_i\subseteq K_{i+1}
	\end{equation}
	for all $i\in \nat$.  For each $i\in \nat$ and each non-Archimedean place $v$ of $K_i$, we seek to define $d_i(v)$ to satisfy the assumptions of Lemma \ref{GeneralExampleConstructor} in such a way that 
	the resulting consistent map $c$ is not $K_i$-Galois invariant for any $i$.  To this end, we first consider the case where there exists a rational prime $q$ such that $x_q\ne 0$.
	
	One easily proves by induction that $K_i/\rat$ is Galois for all $i$ so that we may apply Lemma \ref{NumberFieldChain}.  
	If $v$ is a place of $K_i$ dividing $p$, the first statement of Lemma \ref{NumberFieldChain} ensures that  $\#W_v(K_{i+1}/K_i)$ depends only on $p$.  For simplicity, we shall write
	\begin{equation*}
		M_i(v) = W_v(K_{i+1}/K_i) \quad\mbox{and}\quad M_i = \#W_v(K_{i+1}/K_i).
	\end{equation*}
	By possibly replacing $\{K_i\}$ with a subsequence, we may assume without loss of generality that $M_i > 1$ for all $i\in \nat$.  Clearly the resulting subsequence still satisfies \eqref{FullContain}.
	
	For each prime $p\ne q$, every $i\in \nat$, and every place $v$ of $K_i$ dividing $p$, we let
	\begin{equation} \label{GeneraldDef}
		d_i(v) = x_{p}\frac{[(K_i)_v:\rat_p]}{[K_i:\rat]}.
	\end{equation}
	Therefore, $d_1(p) = x_p$ and 
	\begin{equation} \label{Generald}
		d_i(v) = \sum_{w\in M_i(v)} d_{i+1}(w)
	\end{equation}
	for all $i\in \nat$ and all places $v$ of $K_i$ dividing $p$.
	
	We now proceed by recursively defining $d_i(v)$ for $v$ dividing $q$ starting with $d_1(q) = x_q$.  For each $i\in \nat$ and each place $w$ of $K_{i+1}$ dividing $q$, we assume that
	$w$ divides the place $v$ of $K_i$.  Let $\varepsilon_i(w)\in \rat$ be such that
	\begin{equation} \label{ESum}
		1 = \sum_{w\in M_i(v)} \frac{[(K_{i+1})_w:(K_i)_v]}{[K_{i+1}:K_i]}\varepsilon_i(w)
	\end{equation}
	and
	\begin{equation} \label{EBounds}
		\varepsilon_i(w) \in \left( \exp\left( - 2^{-i} \right),1\right) \cup\left( 1, \exp\left( 2^{-i} \right)\right) 
	\end{equation}
	 for all $w\in M_i(v)$.  Since we have assumed that $M_i >1$, these values for $\varepsilon_i(w)$ necessarily exist, and because of \eqref{EBounds}, they are positive.
	 We now define
	 \begin{equation} \label{SpecialdDef}
	 	d_{i+1}(w) = \frac{[(K_{i+1})_w:(K_i)_v]}{[K_{i+1}:K_i]}\cdot d_i(v)\cdot \varepsilon_i(w),
	\end{equation}
	and by induction on $i$, we quickly see that $d_i(v) \ne 0$ for all $i\in \nat$ and all $v$ dividing $q$.  This definition yields
	\begin{equation*}
		\sum_{w\in M_i(v)} d_{i+1}(w) = d_i(v) \sum_{w\in M_i(v)} \frac{[(K_{i+1})_w:(K_i)_v]}{[K_{i+1}:K_i]}\cdot \varepsilon_i(w),
	\end{equation*}
	and by applying \eqref{ESum}, we obtain
	\begin{equation} \label{Speciald}
	 	d_i(v) = \sum_{w\in M_i(v)} d_{i+1}(w)
	\end{equation}
	for all places $v$ of $K_i$ dividing $q$.  When combining this observation with \eqref{Generald}, the collection of rational numbers $\{d_i(v)\}$ satisfies the assumptions of Lemma \ref{GeneralExampleConstructor}.
	
	We may now assume that $c:\mathcal J\to \rat$ is a consistent map such that $c(K_i,v) = d_i(v)$ for all $i\in \nat$ and all non-Archimedean places $v$ of $K_i$.  If $c$ is $K$-Galois invariant for some number field $K$,
	we may write $K = \rat(\gamma)$ for some $\gamma\in \algt$.  By applying \eqref{FullContain}, there exists $i\in \nat$ such that $K\subseteq K_i$.  Certainly $c$ must also be $K_i$-Galois invariant
	so Lemma \ref{KGaloisInvariant} implies that
	\begin{equation*}
		c(K_{i+1},w) = \frac{[(K_{i+1})_w:(K_i)_v]}{[K_{i+1}:K_i]} \cdot c(K_i,v)
	\end{equation*}
	for all non-Archimedean places $v$ of $K$ and all $w\in M_i(v)$.  However, definition \eqref{SpecialdDef} yields that
	\begin{equation*}
		c(K_{i+1},w) =  \frac{[(K_{i+1})_w:(K_i)_v]}{[K_{i+1}:K_i]} \cdot c(K_i,v) \cdot \varepsilon_i(w).
	\end{equation*}
	Since $c(K_i,v)$ is non-zero, these equalities combine to yield that $\varepsilon_i(w) = 1$, contradicting \eqref{EBounds}.  We have now shown that $c$ cannot be $K$-Galois invariant for any number field $K$.
	Based on our definitions of $d_1(p)$, we have $c(\rat,p) = c(K_1,p) = x_p$ for all rational primes $p$, which completes the proof of the first statement of the Theorem.
	
	To establish the second statement, we shall assume that there exists $B\geq 0$ such that $x_p/\log p \leq B$ for all $p\in M^0_\rat$.   Let $F$ be a number field and $u$ a place of $F$.  
	Then let $i\in \nat$ be such that $F\subseteq K_i$, and thus, the consistency of $c$ implies that
	\begin{equation} \label{BasicCF}
		\left| c(F,u)\right| \leq \sum_{v\in W_u(K_i/F)} \left| d_i(v)\right|.
	\end{equation}
	If $u$ divides a prime $p\ne q$ then definition \eqref{GeneraldDef} yields
	\begin{equation*}
		\left| c(F,u)\right| \leq \sum_{v\in W_u(K_i/F)}|x_p| \frac{[(K_i)_v:\rat_p]}{[K_i:\rat]} = \sum_{v\in W_u(K_i/F)} |x_p|\frac{[(K_i)_v:F_u]}{[K_i:F]} \cdot \frac{[F_u:\rat_p]}{[F:\rat]} 
			= |x_p|\frac{[F_u:\rat_p]}{[F:\rat]},
	\end{equation*}
	and hence
	\begin{equation} \label{BoringBound}
		\left| \frac{c(F,u)}{\log p_u}\right| \leq \frac{|x_p|}{\log p} \cdot \frac{[F_u:\rat_p]}{[F:\rat]}.
	\end{equation}
	
	Next assume that $u$ divides $q$ and let $v$ be an arbitrary place of $K_i$ dividing $u$.  For each $j\leq i$, let $v_j$ be the unique place of $K_j$ such that $v \mid v_j$.
	This means that $v_j \mid v_{j-1}$ for all $2\leq j\leq i$ and that $v_1 = q$.  We now claim that
	\begin{equation} \label{InductionBound}
		\left| d_n(v_n) \right| \leq |x_q|\frac{[(K_n)_{v_n}:\rat_q]}{[K_n:\rat]} \cdot \prod_{j=1}^{n-1} \varepsilon_j(v_{j+1})\quad\mbox{for all } 1\leq n\leq i.
	\end{equation}
	which we shall prove by induction on $n$.  If $n=1$ then both sides are clearly equal to $|x_q|$ so the base case holds.
	Now assume the inductive hypothesis that
	\begin{equation*}
		\left| d_{n-1}(v_{n-1}) \right| \leq |x_q|\frac{[(K_{n-1})_{v_{n-1}}:\rat_q]}{[K_{n-1}:\rat]} \cdot \prod_{j=1}^{n-2} \varepsilon_j(v_{j+1}).
	\end{equation*}
	Using this inequality and \eqref{SpecialdDef} we have
	\begin{align*}
		\left| d_n(v_n) \right|  & = \left| \frac{[(K_n)_{v_n}: (K_{n-1})_{v_{n-1}}]}{[K_n:K_{n-1}]} \cdot d_{n-1}(v_{n-1}) \cdot \varepsilon_{n-1}(v_n)\right| \\
			& \leq  |x_q|\frac{[(K_n)_{v_n}: (K_{n-1})_{v_{n-1}}]}{[K_n:K_{n-1}]} \cdot  \frac{[(K_{n-1})_{v_{n-1}}:\rat_q]}{[K_{n-1}:\rat]} \cdot \left(\prod_{j=1}^{n-2} \varepsilon_j(v_{j+1})\right) \cdot \varepsilon_{n-1}(v_n)\cdot\\
			& = |x_q| \frac{[(K_n)_{v_n}: \rat_q]}{[K_n:\rat]} \cdot \prod_{j=1}^{n-1} \varepsilon_j(v_{j+1})
	\end{align*}
	which establishes \eqref{InductionBound}.  Setting $n=i$ and using the fact that $v = v_i$, we get that
	\begin{equation*}
		\left| d_i(v) \right| \leq |x_p|\frac{[(K_i)_{v}:\rat_p]}{[K_i:\rat]} \cdot \prod_{j=1}^{i-1} \varepsilon_j(v_{j+1}).
	\end{equation*}
	We now apply \eqref{EBounds} to estimate the product on the right hand side of the above inequality.  Specifically, we have
	\begin{equation*}
		\prod_{j=1}^{i-1} \varepsilon_j(v_{j+1}) \leq \prod_{j=1}^{i-1} \exp(2^{-i}) \leq \exp\left( \sum_{i=1}^\infty 2^{-i}\right) = e
	\end{equation*}
	so that
	\begin{equation*}
		\left| d_i(v) \right| \leq e\cdot|x_q|  \frac{[(K_i)_{v}:\rat_q]}{[K_i:\rat]}
	\end{equation*}
	for all places $v$ of $K_i$ dividing $u$.  Combining this inequality with \eqref{BasicCF} we obtain
	\begin{align*}
		\left| c(F,u)\right| & \leq \sum_{v\in W_u(K_i/F)} e\cdot |x_q| \frac{[(K_i)_{v}:\rat_q]}{[K_i:\rat]} \\
			& =  \sum_{v\in W_u(K_i/F)} e\cdot |x_q|\frac{[(K_i)_v:F_u]}{[K_i:F]} \cdot \frac{[F_u:\rat_q]}{[F:\rat]} \\
			& = e\cdot |x_q| \frac{[F_u:\rat_q]}{[F:\rat]},
	\end{align*}
	and therefore
	\begin{equation} \label{InterestingBound}
		\left| \frac{c(F,u)}{\log p_u}\right| \leq \frac{e|x_q|}{\log q} \cdot \frac{[F_u:\rat_q]}{[F:\rat]}.
	\end{equation}
	We have assumed that there exists $B \geq 0$ such that $x_p/\log p \leq B$ for all non-Archimedean places $p$ of $\rat$.  By applying \eqref{BoringBound} and \eqref{InterestingBound}, we have that
	\begin{equation*}
		\left| \frac{c(F,u)}{\log p_u}\right| \leq eB \frac{[F_u:\rat_{p_u}]}{[F:\rat]}
	\end{equation*}
	for all non-Archimedean places $u$ of $F$.  The result follows from Theorem \ref{Continuous} in the case that $x_q \ne 0$ for some rational prime $q$.
	
	Now assume that $x_p = 0$ for all rational primes $p$.  By setting $x'_p = 1$ and applying our previous work, we obtain a consistent map $c':\mathcal J\to \rat$ such that
	\begin{enumerate}[(A)]
		\item $c'(\rat,p) = 1$ for all non-Archimedean places $p$ of $\rat$.
		\item $\Phi_{c'}$ is continuous.
		\item $\stab_G(\Phi_{c'})$ has infinite index in $G$.
	\end{enumerate}
	Now define the map $c'':\mathcal J\to \rat$ by
	\begin{equation*}
		c''(K,v) = -\frac{[K_v:\rat_{p_v}]}{[K:\rat]}.
	\end{equation*}
	Clearly $c''$ is a consistent map with $c(\rat,p) = -1$ for all rational primes $p$, and by Theorem \ref{Continuous}, $\Phi_{c''}$ is continuous.  Moreover, Lemma \ref{KGaloisInvariant} shows that $c''$ is $\rat$-Galois invariant.
	
	Now we define $c = c' + c''$.  We clearly have that $c$ is consistent and that $c(\rat,p) = 0$ for all rational primes $p$.  Additionally, $\Phi_c = \Phi_{c'} + \Phi_{c''}$ so that $\Phi_c$ is continuous.
	Finally, if $\stab_G(\Phi_c)$ has finite index, then $c$ is $K$-Galois invariant for some number field $K$.  Clearly $c''$ is $K$-Galois invariant, and hence, $c' = c - c''$ is $K$-Galois invariant contradicting our assumption
	that $\stab_G(\Phi_{c'})$ has infinite index in $G$.  We have now shown that $c$ satisfies the three required properties completing the proof of the theorem.
\end{proof}

\section{Proofs of Theorems \ref{WeakDense} and \ref{OpenSubgroup}} \label{GaloisApprox}

The proof of Theorem \ref{ClosedSubgroup} required us to find a consistent map $c:\mathcal J\to \rat$ which was not $K$-Galois invariant for any number field $K$.  Of course, our proofs of
Theorems \ref{WeakDense} and \ref{OpenSubgroup} require us to describe the situations in which $c$ is $K$-Galois invariant.  Lemma \ref{KGaloisInvariant} certainly assists with this task, but we also need a strategy to 
construct such maps for any number field $K$.

\begin{lem} \label{InvariantConstructor}
	Suppose that $K$ is a number field, and for each non-Archimedean place $v$ of $K$, let $d_K(v)\in \rat$.  Then there exists a unique $K$-Galois invariant consistent map $c:\mathcal J\to \rat$
	such that $c(K,v) = d_K(v)$ for all non-Archimedean places $v$ of $K$.  Moreover, the following conditions hold:
	\begin{enumerate}[(i)]
		\item\label{GaloisContinuous} $\Phi_c$ is continuous if and only if $\{d_K(v)/\log p_v:v\in M^0_K\}$ is bounded.
		\item\label{CompactSupport} For a set $S$ of rational primes, $c$ is supported on $S$ if and only if $d_K(v) = 0$ for all $v$ not dividing a place in $S$.
	\end{enumerate}
\end{lem}
\begin{proof}
	For each finite extension $L/K$ and each place $w\mid v$, we define
	\begin{equation} \label{dDef}
		d_L(w) = \frac{[L_w:K_v]}{[L:K]} d_K(v)
	\end{equation}
	and note that this definition does indeed agree with the existing values of $d_K(v)$ when $L = K$.  Now assume that $F$ is an arbitrary number field and $u$ is a non-Archimedean place of $F$.
	If $L$ is a number field containing both $K$ and $F$, we let
	\begin{equation} \label{cDef}
		c_L(F,u) = \sum_{w\mid u} d_L(w)
	\end{equation}
	and we claim that $c_L(F,u)$ is independent of $L$.  To see this, assume that $M$ is a number field such that $L\subseteq M$ and adopt the convention that $x$ runs over places of $M$.   Then we have that
	\begin{equation*}
		c_{M}(F,u) = \sum_{x\mid u} d_{M}(x) = \sum_{w\mid u} \sum_{x\mid w} d_M(x).
	\end{equation*}
	For each non-Archimedean place $x$ of $M$, we let $v_x$ be the unique place of $K$ such that $x\mid v_x$.  Then applying definition \eqref{dDef}, we obtain
	\begin{equation*}
		c_{M}(F,u) =\sum_{w\mid u} \sum_{x\mid w}\frac{[M_x:K_{v_x}]}{[M:K]} d_K(v_x).
	\end{equation*}
	For each non-Archimedean place $w$ of $L$, we similarly define $v_w$ to be the unique place of $K$ such that $w\mid v_w$.  Clearly if $x\mid w$ then $v_w = v_x$, and hence, 
	we obtain
	\begin{align*}
		c_{M}(F,u) &  =\sum_{w\mid u} \sum_{x\mid w}\frac{[M_x:K_{v_w}]}{[M:K]} d_K(v_w) \\
			& =\sum_{w\mid u} \sum_{x\mid w}\frac{[M_x:L_w]\cdot [L_w:K_{v_w}]}{[M:L]\cdot [L:K]} d_K(v_w) \\
			& = \sum_{w\mid u} \frac{[L_w:K_{v_w}]}{[L:K]} d_K(v_w) \sum_{x\mid w} \frac{[M_x:L_w]}{[M:L]} \\
			& = \sum_{w\mid u} \frac{[L_w:K_{v_w}]}{[L:K]} d_K(v_w) \\
			& = \sum_{w\mid u} d_L(w) \\
			& = c_L(F,u)
	\end{align*}
	These observations prove that $c_L(F,u)$ is independent of $L$, and therefore, the map $c:\mathcal J\to \rat$ given by $c(F,u) = c_L(F,u)$ is well-defined.  To complete the proof of the first statement, 
	it remains to establish the following four properties of $c$:
	\begin{enumerate}[(A)]
		\item\label{Agreement} $c(K,v) = d_K(v)$ for all non-Archimedean places $v$ of $K$.
		\item\label{GalInv} $c$ is $K$-Galois invariant.
		\item\label{Consistent} $c$ is consistent.
		\item\label{Unique} $c:\mathcal J\to \rat$ is the unique map satisfying \eqref{Agreement}, \eqref{GalInv} and \eqref{Consistent}.
	\end{enumerate}
	
	We note that \eqref{Agreement} follows directly from \eqref{cDef} with $L = F = K$.
	To establish \eqref{GalInv}, we let $v$ a non-Archimedean place of $K$,  $L/K$ be a finite extension, and $w$ a place of $L$ dividing $v$.  Definitions \eqref{dDef} and \eqref{cDef} combine to yield
	\begin{equation} \label{cChain}
		c(L,w) = c_L(L,w) = d_L(w) = \frac{[L_w:K_{v}]}{[L:K]} d_K(v) = \frac{[L_w:K_{v}]}{[L:K]} c_K(K,v) = \frac{[L_w:K_{v}]}{[L:K]}c(K,v),
	\end{equation}
	so Lemma \ref{KGaloisInvariant} shows that $c$ is $K$-Galois invariant.
	
	To prove \eqref{Consistent}, we let $E/F$ be an extension of number fields and assume that $L$ is a number field containing both $E$ and $K$.  As in the proof of Lemma \ref{KGaloisInvariant}, we
	shall assume that $t$ always runs over places of $E$.  If $u$ is a place of $F$ then our definition of $c$ yields
	\begin{equation*}
		c(F,u) = \sum_{w\mid u} d_L(w) = \sum_{t\mid u} \sum_{w\mid t} d_L(w).
	\end{equation*}
	Again applying the definition of $c$, the interior summation is equal to $c(E,t)$ and consistency follows immediately.
	
	Finally, to prove \eqref{Unique}, we assume that $c':\mathcal J\to \rat$ is a $K$-Galois invariant consistent map such that $c'(K,v) = d_K(v)$ for all non-Archimedean places $v$ of $K$.
	Suppose that $F$ is a number field and $L$ is a number field containing both $K$ and $F$.  For each place $u$ of $F$, definitions \eqref{dDef} and \eqref{cDef} yield
	\begin{equation*}
		c(F,u)  = \sum_{w\mid u} d_L(w)  = \sum_{w\mid u} \frac{[L_w:K_{v_w}]}{[L:K]} d_K(v_w)  = \sum_{w\mid u} \frac{[L_w:K_{v_w}]}{[L:K]} c'(K,v_w).
	\end{equation*}
	Since $c'$ is assumed to be $K$-Galois invariant, we have that 
	\begin{equation*}
		c(F,u) = \sum_{w\mid u} c'(L,w),
	\end{equation*}
	and then applying the consistency of $c'$, we conclude that $c(F,u) = c'(F,u)$, as required.
	
	Now we proceed with the proof of \eqref{GaloisContinuous}.  If $\Phi_c$ is continuous then Theorem \ref{Continuous} asserts that there exists $C\geq 0$ such that
	\begin{equation*}
		\left| \frac{c(K,v)}{\log p_v}\right| \leq C\cdot \frac{[K_v:\rat_{p_v}]}{[K:\rat]}.
	\end{equation*}
	As we have that $d_K(v) = c(K,v)$, it follows that $|d_K(v)/\log p_v| \leq C$ for all $v\in M^0_K$.
	
	We now assume that there exists $B\geq 0$ such that $|d_K(v)/\log p_v| \leq B$ for all non-Archimedean places $v$ of $K$.  Suppose that $F$ is a number field and $u$ is a non-Archimedean place of $F$.
	Further assume that $L$ is a number field containing both $F$ and $K$.  Then \eqref{cDef} yields that
	\begin{equation*}
		c(F,u) = \sum_{w\mid u} d_L(w).
	\end{equation*}
	For each place $w$ of $L$, we let $v_w$ be the unique place of $K$ such that $w\mid v_w$.  Now applying \eqref{dDef}, we obtain that
	\begin{equation*}
		c(F,u) = \sum_{w\mid u} \frac{[L_w:K_{v_w}]}{[L:K]} d_K(v_w).
	\end{equation*}
	Therefore, we obtain that
	\begin{align*}
		\left| \frac{c(F,u)}{\log p_u}\right| & \leq B \sum_{w\mid u} \frac{[L_w:K_{v_w}]}{[L:K]}  \\
			& = B\cdot[K:\rat] \sum_{w\mid u} \frac{[L_w:K_{v_w}]}{[L:\rat]} \\
			& \leq B\cdot[K:\rat] \sum_{w\mid u} \frac{[L_w:\rat_{p_u}]}{[L:\rat]} \\
			& = B\cdot[K:\rat] \sum_{w\mid u} \frac{[L_w:F_u]}{[L:F]} \cdot \frac{[F_u:\rat_p]}{[F:\rat]} \\
			& = B\cdot[K:\rat] \cdot \frac{[F_u:\rat_p]}{[F:\rat]},
	\end{align*}
	and the continuity of $\Phi_c$ now follows from Theorem \ref{Continuous}.
	
	Finally, we establish \eqref{CompactSupport}.  If $c$ is supported on $S$ then certainly $d_K(v) = c(K,v) = 0$ for all $v$ not dividing a place in $S$.  Now suppose that $d_K(v) = 0$ for all
	non-Archimedean places $v$ of $K$ not dividing a place in $S$.  Further assume that $F$ is a number field and $u$ is a place of $F$ dividing a place in $S$.  Once again, \eqref{cDef} and \eqref{dDef}
	yield that
	\begin{equation*}
		c(F,u) = \sum_{w\mid u} \frac{[L_w:K_{v_w}]}{[L:K]} d_K(v_w).
	\end{equation*}
	Certainly $v_w$ also divides a place in $S$, so it follows that $c(F,u) = 0$, as required.
\end{proof}

Lemma \ref{InvariantConstructor} provides the fundamental tools that are necessary to establish Theorem \ref{WeakDense}.

\begin{proof}[Proof of Theorem \ref{WeakDense}]
	Suppose that $c\in \mathcal J^*$ is such that $\Phi_c\in \mathcal C^*$.  Therefore, there exists a number field $K$ such that $\Phi_c$ is $K$-Galois invariant.
	Since $c$ is compactly supported, there must exist a finite set $S$ of rational primes such that $c(K,v) = 0$ for all $v$ not dividing a place in $S$.  Therefore, the set $\{c(K,v)/\log p_v:v\in M^0_K\}$ is bounded,
	and it follows from Lemma \ref{InvariantConstructor} that $\Phi_c$ is continuous.

	We now show that $\mathcal C^*$ is countable.
	To this end, we let $\{p_i:i\in \nat\}$ be the complete set of rational primes and let $S_n = \{p_i: 1\leq i\leq n\}$.  For a number field $K$, we let $\mathcal T^*(K,n)$ be the set of all $K$-Galois invariant elements of $\mathcal V^*$
	that are supported on $S_n$.  It is straightforward to check that $\mathcal T^*(K,n)$ is a subspace of $\mathcal V^*$ and that
	\begin{equation} \label{FinComRewrite}
		\mathcal C^* = \bigcup_{K} \bigcup_{n=1}^\infty \mathcal T^*(K,n),
	\end{equation}
	where the outer union runs over all number fields $K$.  As both unions on the right hand side of \eqref{FinComRewrite} are countable, it is enough to show that $\mathcal T^*(K,n)$ is countable.  
	
	Let
	\begin{equation*}
		\mathcal D(K,n) = \left\{ (x_v)_{v\in M^0_K}: x_v\in \rat,\ x_v = 0\mbox{ for all } v\mbox{ not dividing a prime in } S_n\right\}.
	\end{equation*}
	As $S_n$ is a finite set and $\rat$ is countable, $\mathcal D(K,n)$ is also countable.
	Given a point $\xx = (x_v)_{v\in M^0_K} \in \mathcal D(K,n)$, we may apply the first statement of Lemma \ref{InvariantConstructor} with $d_K(v) = x_v$ to obtain a unique $K$-Galois invariant consistent map
	$c_\xx:\mathcal J\to \rat$ such that $c_\xx(K,v) = x_v$ for all non-Archimedean places $v$ of $K$.  According to Lemma \ref{InvariantConstructor}\eqref{CompactSupport}, $c_\xx$ must be supported on $S_n$, and therefore,
	we have a well-defined map $f:\mathcal D(K,n)\to \mathcal T^*(K,n)$ given by 
	\begin{equation*}
		f(\xx) = \Phi_{c_\xx}.
	\end{equation*}	
	
	We now claim that $f$ is surjective.  To see this, let $\Phi\in \mathcal T^*(K,n)$ so that Theorem \ref{AllLinearClassification} yields a unique consistent map $c:\mathcal J\to \rat$ such that $\Phi = \Phi_c$.  
	By definition of $\mathcal T^*(K,n)$, $c$ must also be $K$-Galois invariant.  
	Now let $x_v = c(K,v)$ and $\xx = (x_v)_{v\in M^0_K}$.  By our assumption that $\Phi$ is supported on $S_n$, we know that $x_v = 0$ for all $v$ not dividing a place in $S_n$,
	and hence, $\xx \in \mathcal D(K,n)$.  However, $c_\xx$ is the unique $K$-Galois invariant consistent map such that $c_\xx(K,v) = x_v$, and therefore, $c_\xx = c$.  It follows that $f(\xx) = \Phi_c$ proving that $f$ is surjective
	and that $\mathcal T^*(K,n)$ is countable.
	
	To prove that second statement of the Theorem, we assume that $\Phi\in \mathcal V^*$ and seek a sequence
	$\Phi_n\in \mathcal C^*$ that converges pointwise to $\Phi$.  We let $\{K_i:i\in \nat\}$ be a collection of number fields such that
	\begin{equation*}
		\alg = \bigcup_{i=1}^\infty K_i\quad\mbox{and}\quad K_i \subseteq K_{i+1}
	\end{equation*}
	for all $i\in \nat$.  According to Theorem \ref{AllLinearClassification}, there exists a unique consistent map $c:\mathcal J\to \rat$ such that $\Phi = \Phi_c$.
	For each $n\in \nat$, we define
	\begin{equation} \label{SequentialConsistent}
		d_{K_n}(v) = \begin{cases} c(K_n,v)& \mbox{if } v\mbox{ divides a place in } S_n \\
								0 & \mbox{if } v\mbox{ does not divide a place in } S_n. \end{cases}
	\end{equation}
	According to Lemma \ref{InvariantConstructor}, there exists a unique $K$-Galois invariant consistent map $c_n:\mathcal J\to \rat$ such that $c_n(K_n,v) = d_{K_n}(v)$ for all $v\in M^0_{K_n}$.
	Because of Lemma \ref{InvariantConstructor}\eqref{CompactSupport}, $c_n$ must be compactly supported.  By defining $\Phi_n = \Phi_{c_n}$, we have shown that $\Phi_n\in \mathcal C^*$ for all $n\in \nat$.
	
	Given a point $\alpha\in \algt$, there exists $n\in \nat$ such that $\alpha\in K_n$ and $m\in \nat$ such that $\|\alpha\|_v = 1$ for all
	$v$ not dividing a place in $S_m$.  Letting $N = \max\{m,n\}$ and $i\geq N$ we have
	\begin{equation*}
		\Phi(\alpha) = \Phi_c(\alpha) = \sum_{v} \frac{c(K_i,v)}{\log p_v} \log\|\alpha\|_v = \sum_{v\in S_i}\frac{c(K_i,v)}{\log p_v} \log\|\alpha\|_v.
	\end{equation*}
	Where $v$ runs over places of $K_i$ in each summation.  Definition \eqref{SequentialConsistent} shows that $c(K_i,v) = d_{K_i}(v) = c_i(K_i,v)$ for all $v\in S_i$, and therefore
	\begin{equation*}
		\Phi(\alpha) = \sum_{v\in S_i}\frac{c_i(K_i,v)}{\log p_v} \log\|\alpha\|_v = \sum_{v}\frac{c_i(K_i,v)}{\log p_v} \log\|\alpha\|_v = \Phi_i(\alpha).
	\end{equation*}
	We have now shown that $\Phi(\alpha) = \Phi_i(\alpha)$ for all $i\geq N$ completing the proof.
\end{proof}

The proof of Theorem \ref{OpenSubgroup} is inspired by an example how Lemma \ref{InvariantConstructor} can be used to build a continuous $K$-Galois invariant extension of $\Omega$ that is distinct from \eqref{CanonicalExtension}.  
Let $q(x) = x^2 + 1$ and note that by Hensel's Lemma (see \cite{Koch}, for example), $q$ has a root $i\in \rat_5$.  Now let $K = \rat(i)$ and write $s_v = [K_v:\rat_{p_v}]/[K:\rat]$ so that
\begin{equation*}
	1 = \sum_{v \mid p} s_v
\end{equation*}
for all places $p$ of $\rat$.  For each place $v$ of $K$ dividing $5$, we have that $s_v = 1/2$.  By applying \eqref{DegreeSum}, there must exist exactly two places $v_1$ and $v_2$ of $K$ dividing $5$.   
We define
\begin{equation} \label{AltdValues}
	d_K(v) = \begin{cases} -s_v & \mbox{if } v\not\in\{v_1,v_2\} \\
			-1/3 & \mbox{if } v = v_1 \\
			-2/3 & \mbox{if } v = v_2 \end{cases}
\end{equation}
and let $c:\mathcal J\to \rat$ be the $K$-Galois invariant consistent map described by Lemma \ref{InvariantConstructor}.  Clearly $c(\rat,p) = -1$ for all primes $p$, and therefore, Theorem \ref{OmegaExtension}
shows that $\Phi_c$ really is an extension of $\Omega$ while Lemma \ref{InvariantConstructor}\eqref{GaloisContinuous} shows it to be continuous.  However,  we also have that
\begin{equation*}
	c(K,v_1) \ne \frac{[K_{v_1}:\rat_5]}{[K:\rat]} c(\rat,p)
\end{equation*}
so according to Lemma \ref{KGaloisInvariant}, $c$ cannot be $\rat$-Galois invariant.  It now follows that $\stab_G(\Phi_c) = \gal(\alg/\rat(i))$ and that $\Phi_c$ is distinct from \eqref{CanonicalExtension}.

\begin{proof}[Proof of Theorem \ref{OpenSubgroup}]
	Suppose that $K$ is a number field such that $\gal(\alg/K) = H$.    We shall now select points $d_K(v)\in \rat$ for use in Lemma \ref{InvariantConstructor}.
	If $K = \rat$ then we let $d_\rat(p) = x_p$ for all $p\in M^0_\rat$.  Now Lemma \ref{InvariantConstructor} gives a $\rat$-Galois invariant consistent map $c:\mathcal J\to \rat$ 
	such that $c(\rat,p) = x_p$ for all $p\in M^0_\rat$.   Vacuously,  $c$ cannot be $F$-Galois invariant for any proper subextension $F/\rat$.  If $\{x_p/\log p:p\in M^0_\rat\}$ then 
	Lemma \ref{InvariantConstructor}\eqref{GaloisContinuous} ensures that $\Phi_c$ is continuous,
	
	Now assume that $K\ne \rat$ and let $K_1,K_2,\ldots,K_n$ be the complete list of maximal subfields of $K$.
	For each such $i$, the Chebotarev density theorem \cite{FriedJarden} implies the existence of infinitely many places of $K_i$ 
	that split completely in $K$.  Therefore, for each $1\leq i\leq n$, we may assume that $v_i$ is a non-Archimedean place of $K_i$ satisfying the following two properties:
	\begin{enumerate}[(A)]
		\item\label{SplitsCompletely} $v_i$ splits completely in $K$
		\item\label{NonCompliance} if $i\ne j$ then $v_i$ and $v_j$ divide distinct places of $\rat$.
	\end{enumerate}
	For every $1\leq i\leq n$, property \eqref{SplitsCompletely} means that there are at least two places of $K$ dividing $v_i$.  For each $i$ and each place $v\mid v_i$, 
	we may select $\varepsilon_i(v) \in \rat\setminus\{0\}$ such that
	\begin{equation} \label{Balance}
		\sum_{v\mid v_i} \varepsilon_i(v) = 0.
	\end{equation}
	Because of property \eqref{NonCompliance},  each place $v$ of $K$ divides $v_i$ for at most one value of $i$, and therefore, $\varepsilon_i(v)$ depends only on $v$.  By defining
	\begin{equation*}
		\widehat{\mathcal M} = \left\{ v\in M^0_K:  v\mid v_i\mbox{ for some } 1\leq i\leq n\right\}
	\end{equation*}
	we may now let $\varepsilon(v) = \varepsilon_i(v)$ for all places $v\in \widehat{\mathcal M}$.  If $v$ is a non-Archimedean place of $K$, we define
	\begin{equation} \label{DefinitionOfd}
		d_K(v) = x_{p_v}\frac{[K_v:\rat_p]}{[K:\rat]} + \begin{cases} \varepsilon(v) & \mbox{if } v\in\widehat{\mathcal M} \\ 0  & \mbox{if } v\not\in \widehat{\mathcal M}.\end{cases}
	\end{equation}
	Clearly $d_K(v)\in \rat$ for all $v\in M^0_K$, so Lemma \ref{InvariantConstructor} yields a $K$-Galois invariant consistent map $c:\mathcal J\to \rat$ such that $c(K,v) = d_K(v)$ for all $v\in M^0_K$.
	We claim that $c$ satisfies the required properties.
	
	To prove \eqref{RationalFixed}, we let $p$ be a non-Archimedean place of $\rat$.   If there does not exist $1\leq i\leq n$ such that $v_i\mid p$ then \eqref{DefinitionOfd}, then $v$ cannot belong to $\widehat{\mathcal M}$ for
	any place $v\mid p$.  Now using the consistency of $c$, we deduce that
	\begin{equation*}
		c(\rat,p) = \sum_{v\mid p} c(K,v) = x_p\sum_{v\mid p} \frac{[K_v:\rat_p]}{[K:\rat]} = x_p,
	\end{equation*}
	establishing \eqref{RationalFixed} in this case.  Alternatively, if $1\leq i\leq n$ is such that $v_i\mid p$, then property \eqref{NonCompliance} means that $v_j$ cannot divide $p$ for any $j\ne i$.
	Hence, we deduce from \eqref{DefinitionOfd} that
	\begin{equation} \label{cOmega}
		c(\rat,p) = \sum_{v\mid p}c(K,v) = \sum_{v\mid v_i} \left(x_p\frac{[K_v:\rat_p]}{[K:\rat]} + \varepsilon(v)\right) + \sum_{\substack{ v\mid p \\ v\nmid v_i}}x_p\frac{[K_v:\rat_p]}{[K:\rat]}
	\end{equation}
	Now using \eqref{Balance} we again obtain that $c(\rat,p) = x_p$, establishing \eqref{RationalFixed}.
	
	We now prove \eqref{StabilizerFixed} by contradiction, so we assume that $1\leq i\leq n$ is such that $c$ is $K_i$-Galois invariant.   Now using the consistency of $c$ along with \eqref{Balance} and \eqref{DefinitionOfd} we obtain that
	\begin{align*}
		c(K_i,v_i) & = \sum_{v\mid v_i} c(K,v) \\
			& = \sum_{v\mid v_i}\left(  x_{p_v}\frac{[K_v:\rat_{p_v}]}{[K:\rat]} + \varepsilon(v)\right) \\
			& = x_{p_v} \sum_{v\mid v_i} \frac{[K_v:(K_i)_{v_i}]}{[K:K_i]}\cdot \frac{[(K_i)_{v_i}:\rat_{p_v}]}{[K_i:\rat]} \\
			& = x_{p_v} \frac{[(K_i)_{v_i}:\rat_{p_v}]}{[K_i:\rat]}
	\end{align*}
	Now fix a place $v$ of $K$ dividing $v_i$.  Since we have assumed that $c$ is $K_i$-Galois invariant, Lemma \ref{KGaloisInvariant} yields that
	\begin{equation*} \label{cGalois}
		c(K,v) = \frac{[K_v:(K_i)_{v_i}]}{[K:K_i]} c(K_i,v_i),
	\end{equation*}
	so we obtain that
	\begin{equation*}
		c(K,v) = x_{p_v} \frac{[K_v:(K_i)_{v_i}]}{[K:K_i]} \cdot \frac{[(K_i)_{v_i}:\rat_{p_v}]}{[K_i:\rat]} = x_{p_v} \frac{[K_v:\rat_{p_v}]}{[K:\rat]}.
	\end{equation*}
	However, \eqref{DefinitionOfd} gives
	\begin{equation*}
		c(K,v) =x_{p_v} \frac{[K_v:\rat_{p_v}]}{[K:\rat]} + \varepsilon(v),
	\end{equation*}
	contradicting our assumption that $\varepsilon(v) \ne 0$.
	
	Finally, we assume that $\{x_{p_v}/\log p_v:v\in M^0_K\}$ is bounded.  We certainly have that $|d_K(v)/\log p_v| \leq |x_{p_v}/\log p_v|$ for all $v\not\in \widehat{\mathcal M}$, and since $\widehat{\mathcal M}$
	is a finite set, we must have that $\{d_K(v)/\log p_v:v\in M^0_K\}$ is bounded.  Now Lemma \ref{InvariantConstructor}\eqref{GaloisContinuous} shows that $\Phi_c$ is continuous.
\end{proof}

\section{Acknowledgment}

The author wishes to thank Professor J.D. Vaaler for his many useful comments regarding this work including the recommendation to pursue the problem of classifying $\rat$-valued linear functionals on
$\mathcal V$.  Additionally, we thank the referee for observing the alternate proof of Theorem \ref{AllLinearClassification} that involves inverse limits as well as many additional helpful remarks.


\providecommand{\bysame}{\leavevmode\hbox to3em{\hrulefill}\thinspace}
\providecommand{\MR}{\relax\ifhmode\unskip\space\fi MR }
\providecommand{\MRhref}[2]{%
  \href{http://www.ams.org/mathscinet-getitem?mr=#1}{#2}
}
\providecommand{\href}[2]{#2}

\end{document}